\documentclass[12pt,leqno,a4paper]{amsart}
\usepackage{amssymb,enumerate}
\newcommand{\Irr}{{\operatorname{Irr}}}

\newcommand{\Sz}{\operatorname{Sz}}
\newcommand{\PSp}{\operatorname{PSp}}
\newcommand{\tr}{\operatorname{tr}}
\renewcommand{\phi}{\varphi}

\theoremstyle{definition}

\theoremstyle{remark}

\begin{document}

\title[A counterexample to Feit's Problem VIII on decomposition numbers]{A counterexample to Feit's Problem VIII
on decomposition numbers}

\date{\today}
 
\author{Gabriel Navarro}
\address{Department of Mathematics, Universitat de Val\`encia, 46100 Burjassot,
        Spain.}
\email{gabriel@uv.es}
\author{Benjamin Sambale}
\address{FB Mathematik, TU Kaiserslautern, Postfach 3049,
        67653 Kaisers\-lautern, Germany.}
\email{sambale@mathematik.uni-kl.de}

\begin{abstract}
We find a counterexample to Feit's Problem VIII on the bound of decomposition numbers. This also answers a question raised by 
T. Holm and W. Willems.
\end{abstract}

\thanks{ The first author is partially supported by the
Spanish Ministerio de Educaci\'on y Ciencia Proyectos  MTM2016-76196-P  and
Prometeo II/Generalitat Valenciana. He also thanks Gerhard Hiss for useful discussions
on this, while at the ICTS in Bangalore, and the ICTS center for the hospitality.
The second author thanks the German Research
Foundation SA 2864/1-1 and Daimler Benz Foundation 32-08/13.
Both authors thank Thomas Breuer and Klaus Lux for conversations on this subject.}



\maketitle

Richard Brauer  asked if the Cartan invariants $c_{\phi\psi}$ of a $p$-block $B$ of a finite group $G$ are at most $p^d$ where $d$ is the defect of $B$. It is well-known that 
Peter Landrock showed that the Suzuki group $\Sz(8)$ with $p=2$ is a counterexample to Brauer's question. 
Since 
\[c_{\phi\phi}=\sum_{\chi\in\Irr(B)}{(d_{\chi\phi})^2},\]
Walter Feit, in his list of open problems in Representation Theory,
 proposed the following weaker question on decomposition numbers in his book \cite[Problem (VIII), p.\,169]{Feit}:

\medskip

(VIII)~~{\it Is $(d_{\chi \varphi})^2 \le p^d$ whenever $\chi, \varphi$ lie in a 
block of defect $d$?}

\medskip

To the present authors' surprise apparently no one has noticed that the Atlas of Brauer characters~\cite{BrauerAtlas} contains two counterexamples to Feit's problem: $\PSp_4(4).4$ and $\Sz(32).5$ both for $p=2$, both in the principal block. In the first case, $44$ occurs as a decomposition number and in the second case, $47$ occurs (see \cite{MOC} for instance). For both groups we have $|G|_2=2^{10}$.
 
It is interesting to speculate on whether or not (VIII) (or even Brauer's original question) has a positive answer for odd primes.
This would prove, using Brauer's $k(B)$-conjecture, that 
$c_{\phi\psi}\le\max\{c_{\phi\phi},c_{\psi\psi}\}\le p^{2d}$. In fact, we are not aware of any counterexample
to this, even for $p=2$. Also, we are not aware of any example where
the inequality $d_{\chi \phi} \le p^d$ does not hold.

\bigskip

Another way of relaxing Brauer's question was proposed by Holm and Willems. They ask in \cite[p.\,594]{HW} if
\[\tr C\le l(B)p^d\] 
holds for every $p$-block $B$ with defect $d$ and Cartan matrix $C$. An easy computation with GAP~\cite{GAP} shows that our examples above provide counterexamples also for this question.


  

\end{document}